\documentclass[12pt ]{article}
\title{Jacobi fields and odular structure of affine manifolds}
\author{Alexander I. Nesterov
\thanks{Departamento de F\'\i sica, Universidad de Guadalajara, Guadalajara,
Jalisco, M\'exico;  E-mail: nesterov@cencar.udg.mx.}}

\date{\today}

\usepackage{amsmath,amsthm}

\chardef\bslash=`\\ % p. 424, TeXbook
\newtheorem{thm}{Theorem}[section]

\newtheorem{lem}[thm]{Lemma}
\newtheorem{prop}[thm]{Proposition}

\theoremstyle{definition}
\newtheorem{defn}{Definition}[section]

\theoremstyle{remark}
\newtheorem{rem}{Remark}[section]

\newcommand{\A}{\mathcal{A}}
\newcommand{\B}{\mathcal{B}}
\newcommand{\eval}[2][\right]{\relax
  \ifx#1\right\relax \left.\fi#2#1\rvert}

\newcommand{\Bvert}[1]{\left.#1\right\rvert}

\input amssym.def
\input amssym

\begin{document}
\date{~}
\maketitle

\begin{abstract}
The connection between Jacobi fields and odular structures of affine
manifold is established. It is shown that the Jacobi fields generate the
natural geoodular structure of affinely connected manifolds.
\end{abstract}



%\markboth{Sample paper for the {\protect\ntt\lowercase{amsmath}} package}
%{Sample paper for the {\protect\ntt\lowercase{amsmath}} package}
\renewcommand{\sectionmark}[1]{}

\newpage
\section{Introduction}

The recent development of geometry has shown the importance
of non-associative algebraic structures, such as quasigroups, loops and
odules. For instance, it is possible to say that the nonassociativity is
the algebraic equivalent of the differential geometric concept of the
curvature. The corresponding construction may be described as follows. In
a neighbourhood of an arbitrary point on a manifold with the affine
connection one can introduce the geodesic local loop, which is uniquely
defined by means of the parallel translation of
geodesics along geodesics (Kikkawa, 64; Sabinin 72, 77). The family of local
loops constructed in this way uniqeuly defines the space with affine
connection, but not every family of geodesic loops on a manifold yields an
affine connection. It is necessary to add some algebraic identities
connecting loops in different points. This additional algebraic  structures
(so-called geoodular structures) were introduced  and the
equivalence of the categories of geoodular structures and of affine
connections was shown by Sabinin (77, 81).  The development of this
line gave a start to a new approach to manifold with affine coonection -
{\it Loopuscular and Odular Geometry.} The main algebraic structures
arising in this approach are related to nonassociative algebra and
theory of quasigroups and loops.

In our paper we investigate the connection between Jacobi fields and odular
structure of affinely connected manifold. In Section 2 we outline the main
consructions of the loopuscular and odular geometry. In Section 3 we proof
our main result ({\bf Theorem 3.3}).

\section{Odular structure of manifolds with affine connection}
Below we outline the main facts from the {loopuscular geometry} following
(Sabinin, 77, 81, 87, 88, 98).

\defn{Let $\langle{\frak M},\cdot,\varepsilon\rangle$ be a partial {\it
magma} with a binary operation $(x,y)\mapsto x\cdot y$ and the
neutral element $\varepsilon,\; x\cdot\varepsilon =\varepsilon\cdot
x =x$; $\frak M$ be a smooth manifold (at least $C^1$-smooth) and
the operation of multiplication (at least $C^1$-smooth) be defined
in some neighbourhood $U_\varepsilon$,} then  $\langle {\frak
M},\cdot,\varepsilon\rangle$ is called a {\it partial loop on
$\frak M$.}

\rem{The operation of multiplication is locally left and right
invertible. This means, if $x\cdot y = L_x y= R_y x$, then there exist
$L^{-1}_x$ and $R^{-1}_x$ in some neighbourhood of the neutral element
$\varepsilon$:
\[
L_a(L^{-1}_a x)=x, \quad R_a(R^{-1}_a x)=x.
\]
}

{The vector fields $A_j$ defined in $U_\varepsilon$ by
\begin{equation}
A_j(x)=\bigl((L_x)_{*,\varepsilon}\bigr)^i_j\frac{\partial}{\partial
x^i}
\end{equation}
are called the {\it left basic fundamental fields}. Similarly, the
{\it right basic fundamental fields} $B_j$ are defined by
\begin{equation}
B_j(x)=\bigl((R_x)_{*,\varepsilon}\bigr)^i_j\frac{\partial}{\partial
x^i}.
\end{equation}
}

The solution of the equation
\begin{equation}
\frac{df^i(t)}{dt} =L^i_j(f(t))x^j, \quad f(0) = \varepsilon
\end{equation}
yields the {\it exponential map}
\[
{\rm Exp}: X\in T_\varepsilon ({\frak M})\longrightarrow
{\rm Exp}X \in{\frak M}.
\]

The unary operation
\begin{equation}
tx={\rm Exp}t {\rm Exp}^{-1}x,
\end{equation}
based on the exponential map, is called the {\it left canonical
unary operation} for $\langle {\frak M},\cdot,\varepsilon\rangle$.
A smooth loop $\langle {\frak M},\cdot,\varepsilon\rangle$
equipped with its canonical left unary operations is called the
{\it left canonical preodule} $\langle {\frak
M},\cdot,(t)_{t\in\Bbb R}, \varepsilon\rangle$. If one more operation is
introduced
\begin{equation} x+y ={\rm Exp}({\rm Exp}^{-1}x + {\rm
Exp}^{-1}y),
\end{equation}
then we obtain the {\it canonical left prediodule} of a loop, $\langle
{\frak M},\cdot,+,(t)_{t\in\Bbb R}, \varepsilon\rangle $. {A canonical left
preodule (prediodule) is called the {\it left odule (diodule)}, if the {\it
monoassociativity} property
\begin{equation}
tx\cdot ux = (t+u)x
\end{equation}
is satisfied.} In the smooth case for an odule the left and the right
canonical operations as well as the exponential maps coincide.

\defn{Let $\frak M$} be a smooth manifold and
\[
L:  (x,y,z)\in {\frak M} \mapsto L(x,y,z)\in {\frak M}
\]
a smooth partial ternary operation, such that $x_{\dot a}y =L(x,a,z)$
defines  in the some neighbourhood of the point $a$   the loop with the
neutral $a$, then the pair $\langle {\frak M}, L\rangle$  is called a {\it
loopuscular structure (manifold)}.

{A smooth manifold  $\frak M$ with a smooth partial ternary operation
$L$ and smooth binary operations $\omega_t : (a,b)\in {\frak M}\times{\frak
M} \mapsto \omega_t(a,b)=t_a b\in {\frak M}, \;(t\in \Bbb R)$, such that
$x_{\dot a}y =L(x,a,y)$ and $t_a z = \omega_t(a,z)$ determine in some
neighbourhood of an arbitrary point $a$ the odule with the neutral element
$a$, is called a {\it left odular structure (manifold)} $\langle {\frak M},
L, (\omega_t)_{t\in {\Bbb R}}\rangle$. Let $\langle {\frak M}, L,
(\omega_t)_{t\in {\Bbb R}}\rangle$ and $\langle {\frak M}, \Lambda,
(\omega_t)_{t\in {\Bbb R}}\rangle$ be odular structures then $\langle
{\frak M},L, \Lambda, (\omega_t)_{t\in {\Bbb R}}\rangle$ } is called a {\it
diodular structure (manifold). If $x_{\stackrel{+}a}y = \Lambda(x,a,y)$}
and $t_a x =\omega_t(a,x)$ define a vector space, then such a diodular
structure is called a {\it linear diodular structure}

\defn{\rm(Sabinin, 77, 81)} {Let $\frak M$ be a $C^k$-smooth $(k\geq 3)$
affinely connected manifold and the following operations are defined  on
$\frak M$:
\begin{align}
& L(x,a,y)= {x\; {}_{\stackrel{\cdot}{a}}\; y}
={\rm Exp}_x \tau^a_x {\rm Exp}^{-1}_a y, \\
& \omega_t(a,z) = t_a z= {\rm Exp}_a t {\rm Exp}^{-1}_a z, \\
&\Lambda (x,a,y) = {x\; {}_{\stackrel{+}{a^{}}}\; y}
={\rm Exp}_a ({\rm Exp}^{-1}_a x + {\rm Exp}^{-1}_a y),
\end{align}
${\rm Exp}_x$ being the exponential map at the point $x$ and $\tau^a_x$
the parrallel translation along the geodesic going from $a$ to $x$.
The construction above is called a {\em natural linear geodiodular
structure of an affinely connected manifold} $({\frak M},
\nabla)$.}

\rem {Any $C^k$-smooth $(k\geq 3)$  affinely connected manifold can be
considered as a geoodular structure.}

\defn{(\rm Sabinin, 77, 81, 86)} {Let $\langle {\frak M}, L\rangle$ be a
loopuscular structure of a smooth manifold $\frak M$. Then the formula
\begin{gather*}
\nabla_{X_a}Y=\left\{\frac{d}{dt}
\left([(L^a_{g(t)})_{*,a}]^{-1}Y_{g(t)}\right)\right\}_{t=0},\\
g(0)=a, \quad \dot g(0)=X_a,
\end{gather*}
$Y$ being a vector field in the neighbourhood of a point $a$, defines the
{\em tangent affine connection}.}

In coordinates the components of affine connection are written as
\[
\Gamma^i_{jk}(a) = -\left[\frac{\partial^2 (x{}_{\dot a}{}y)^i}{\partial x^j
\partial x^k}\right]_{x=y=a}
\]

The equivalence of the categories of geoodular (geodiodular) structures and
of affine connections is formulated as follows (Sabinin, 77, 81).

{\prop The tangent affine connection $\overline \nabla$ to
the natural geoodular (geodiodular) structure of an affinely connected
manifold $({\frak M}, \nabla)$  coincides with $\nabla$.}

{\prop The natural geoodular (geodiodular) structure
$\overline {\frak M}$ of the tangent affine connection to a natural geoodular
(geodiodular) structure $\frak M$  coincides with $\frak M$.}

\section{Jacobi fields and odular structure}

A vector field $X$ along a geodesic $\gamma$ is called a {\it Jacobi field}
if it satisfies the Jacobi differential equation (the geodesic deviation
equation)
\begin{equation}
\frac{D^2 X}{dt^2} + {\mbox{\boldmath $R$}}(X,Y)Y=0,
\end{equation}
$D/dt$ being the operator of the parallel translation along $\gamma$,
${\mbox{\boldmath $R$}}(X,Y)$ the curvature operator and
$Y=d\gamma/dt$ the tangent vector to the geodesic $\gamma(t)$.

As is known the Jacobi equation has $2n$ linearly independent
solutions, which are completely determined by the initial
conditions: $X(0)$ and \linebreak $ D{X(0)}/{dt}\in
T_{\gamma}({\frak M})$. Besides, every geodesic admits two natural
Jacobi fields. The first one is defined as follows: $X_1=Y$, and
the second one: $X_2=tY$; $t$ being canonical parameter along the
geodesic (Kobayashi and Nomizu, 69; Milnor, 73).

\defn{\rm (Kobayashi and Nomizu, 69)} {A one-parametric family of geodesics
$\alpha(s,t), \; -\varepsilon < s < \varepsilon$, such that
$\alpha(0,t)=\gamma(t), \; 0\leq t \leq 1$, is called a {\it variation} of
the geodesic} $\gamma(t))$.

This means, that there is a $C^{\infty}$-differentiable map from
$[0,1]\times (-\varepsilon,\varepsilon)$ to $\frak M$ such that:

(i) for each given $s\in (-\varepsilon, \varepsilon) \quad
\alpha(s,t)$ is the geodesic;

(ii) $\alpha(0,t) = \gamma(t)$ for $0\leq t\leq 1$.

An {\it infinitesimal variation} $X$ of the geodesic $\gamma(t)$ is defined
as follows:
\begin{eqnarray}
\Bvert{X}_{x} = \frac{\partial\alpha(s,t)}{\partial s} \quad
\mbox{for} \quad 0\leq t \leq 1,
\end{eqnarray}
where $\alpha(s,t)$ is a variation. Further the following theorem is
important.

{\thm {\rm (Kobayashi and Nomizu 69)} A vector field $X$ along the
geodesic $\gamma$ is a Jacobi field if and only if it is an
infinitesimal variation for $\gamma$.}

Let us consider a geodesic variation $\alpha( s,t)$ of a geodesic
$\gamma(t)$ starting at a point $a$ such, that
$\Gamma( s)=\alpha( s,0)$ be a geodesic passing through the point $a$, as
well, and the infinitesimal variation  $X=\partial\alpha( s,t)/\partial  s$
satisfies
\begin{eqnarray}
\Bvert{X}_{\Gamma}\equiv X( s,0) \neq 0,\quad
\frac{D X( s,0)}{dt}\equiv \Bvert{\frac{D
X}{dt}}_{\Gamma}=0.
\end{eqnarray}
Let $\zeta =\partial\alpha( s,t)/\partial  s |_{ s=t=0}$ and $\xi
=\partial\alpha( s,t)/\partial t |_{ s=t=0}$  be tangent vectors at
the point $a$ to geodesics $\Gamma( s)$ and $\gamma(t)$ respectively.
The following lemma is valid.

{\lem Let $\tau^a_x $ be the parallel translation along $\Gamma$, then it
holds
\[
\alpha( s,t)={\rm Exp}_{x( s)} \tau^a_{x( s)} {\rm Exp}^{-1}_a y(t),
\]
where $x(s)={\rm Exp}_a( s\zeta)$ and $y(t)={\rm Exp}_a(t\xi)$.\label{Lem}
}

\begin{proof}
Any geodesic $\tilde \alpha(x,t)$ pasing through a point $x\in
\Gamma$ can be presented as $\tilde\alpha(x,t)= {\rm Exp}_x (t \eta)$,
$\eta$ being the tangent vector to $\tilde\alpha$ at the point $x$,
and $t$ the canonical parameter. Representing $\eta$ as  $\eta = \tau^a_x
\xi$, where $\tau^a_x$ is the parallel translation along $\Gamma$, we find
\[
\tilde\alpha(x,t)= {\rm Exp}_x (t \tau^a_x \xi ) ={\rm
Exp}_x \tau^a_x (t\xi ) .
\]
Introducing $y(t) = {\rm Exp}_a(t\xi)$ and $x=x(s)=\alpha(s,0)$, we obtain
\[
\tilde\alpha(x(s),t)= {\rm Exp}_{x(s)} (t \tau^a_{x(s)} \xi )
={\rm Exp}_{x(s)} \tau^a_{x(s)} {\rm Exp}^{-1}_a y(t).
\]
On the other hand, $\tilde\alpha(x(s),t)$ can be obtained as the result of
the variation of the geodesic $\gamma(t)$. This yields
$\tilde\alpha(x(s),t)=\alpha(s,t)$ and, obviously,
\[
\alpha( s,t)={\rm Exp}_{x( s)} \tau^a_{x( s)} {\rm Exp}^{-1}_a y(t).
\]
The proof is completed.
\end{proof}

{\thm {\rm (Nesterov, 1989)} Jacobi fields generate the natural linear
geodiodular structure of an affinely connected space $({\frak M},\nabla)$.
}

\begin{proof}

Taking into account Lemma \ref{Lem}, we define the operation
$L(x(s),a,y(t))$ in the following way:
\[
L(x(s), a,y(t))= \alpha(s,t)=
={\rm Exp}_{x(s)} \tau^a_{x(s)} {\rm Exp}^{-1}_a y(t).
\]
Actually, this operation is defined for arbitrary points $x$ and $y$ in
the some neighbourhood $U_a$ of the point $a$. Thus, it can be
represented as
\[
L(x, a,y)= {\rm Exp}_{x} \tau^a_{x} {\rm Exp}^{-1}_a y,
\]
(see {\bf Def. 2.3}, (7)).

Now let us consider the following variation of geodesic $\gamma(t)$:
\[
\beta(s,t)={\rm Exp}_a(t\xi +st\eta).
\]
The infinitesimal variation  $X=\partial\beta( s,t)/\partial s$
satisfies
\[
X( s,0) = 0,\quad \frac{D X( s,0)}{dt}\neq 0.
\]
Defining $x(t)={\rm Exp}_a(t\xi)$ and $y(s,t)={\rm Exp}_a(st\eta)$,
we obtain
\begin{equation}
\beta(s,t) = {\rm Exp}_a({\rm Exp}^{-1}_a x(t)
+ {\rm Exp}^{-1}_a y(s,t)).
\end{equation}
In fact, this is valid for arbitrary points $x,y \in U_a$. Thus, one
can identify the operation (9) as the geodesic variation $\beta$:
\[
\Lambda (x,a,y) = \beta(s,t)
={\rm Exp}_a ({\rm Exp}^{-1}_a x + {\rm Exp}^{-1}_a y),
\]

The operation (8) is related to the existence of the affine parameter
along a geodesic. The proof is completed.

\end{proof}
\section*{Acknowledgments}

The author is grateful to L.V. Sabinin for helpful discussions and
comments. He is also greatly indebted to L.L. Sbitneva for reading the
manuscript and helpfull remarks. This work was supported partly by CONACyT,
Grant No.  1626P-E.

\section*{References}
\begin{description}
\item[] {} Kikkawa, M.: {\it On local loops in affine manifolds.} J.
Sci.  Hiroshima Univ. Ser A-I Math. {\bf 28}, 199 (1961).

\item[] {} Kobayashi, S. K. and Nomizy, K.:
{\it Foundations of Differential Geometry}, Vol. 1. New York:
Interscience Publisher 1963.

\item[] {} Kobayashi, S. K. and Nomizy, K.:
{\it Foundations of Differential Geometry}, Vol. 2. New York: Interscience
Publisher 1969.

\item[]{} Milnor, J.: {\it Morse theory.} Princeton: Princeton
Univ. Press 1973.

\item[] {} Nesterov, A. I.: {\it Methods of nonassociative algebra in
physics.} Doctor of Science Theses. Tartu: 1989.

\item {}
Sabinin, L. V.: {\it  On the equivalence of categories of loops and
homogeneous spaces}, Soviet Math. Dokl. {\bf 13}, 970 (1972a).

\item[] {}Sabinin, L. V.: {\it  The geometry of loops},
 Mathematical Notes, {\bf 12}, 799 (1972b).

\item[] {}Sabinin, L. V.:
{\it  Odules as a new approach to a geometry with a connection,}
Soviet Math. Dokl. {\bf 18}, 515 (1977).

\item[] {}Sabinin, L. V.:
{\it Methods of Nonassociative Algebra in Differantial Geometry},
in {Supplement to Russian translation of S.K. Kobayashi and K. Nomizy
``Foundations of Differential Geometry'', Vol. 1.} Moscow: Nauka 1981.

\item[] {}Sabinin, L. V.:
{\it Differential equations of smooth loops}, in:
 {Proc. of Sem. on Vector and Tensor Analysis,}
{\bf 23} 133. Moscow: Moscow Univ. 1988.

\item[]{} Sabinin, L. V.: {\it  Differential Geometry and Quasigroups},
Proc. Inst. Math. Siberian Branch of Ac. Sci. USSR, {\bf 14}, 208 (1989).

\item[] {}Sabinin, L. V.: {\it On differential equations of smooth loops},
Russian Mathematical Survey, {\bf 49} 172 (1994).

\item[] {}Sabinin, L. V.:
{\it Smooth quasigroups and loops.} Dordrecht: Kluwer Academic
Publishers 1998.

\end{description}

\end{document}